\documentclass[reqno, 11pt]{amsart}

\usepackage[top=4cm, bottom=3cm, left=3.2cm, right=3.2cm]{geometry}
\newtheorem{thm}{Theorem}

\newtheorem{lem}{Lemma}
\newtheorem{prop}{Proposition}
\theoremstyle{definition}

\theoremstyle{remark}

\newcommand{\dint}{\displaystyle\int}
\newcommand{\dlim}{\displaystyle\lim}
\newcommand{\dsup}{\displaystyle\sup}

\newcommand{\dsum}{\displaystyle\sum}

\newcommand{\F}{\mathcal{F}}
\newcommand{\R}{\mathbb{R}}
\newcommand{\T}{\mathbb{T}}
\newcommand{\I}{\mathbb{I}}

\numberwithin{equation}{section} \numberwithin{lem}{section}
\numberwithin{thm}{section} \numberwithin{prop}{section}
\numberwithin{cor}{section} \numberwithin{rem}{section}

\begin{document}
\title[Incompressible limits for charge-carrier
transport]{Incompressible type limit analysis of a hydrodynamic model for charge-carrier
transport}
%\tableofcontents
\author{Li Chen}
\address{Department of Mathematical Sciences\\
 Tsinghua University\\
  Beijing, 100084, PeopleÕs Republic of China}
  \email{lchen@math.tsinghua.edu.cn}
\author{Donatella Donatelli}
\address{Dipartimento di Matematica Pura and Applicata\\
Universit\`a di L'Aquila\\
67100 L'Aquila, Italy}
\email{donatell@univaq.it}
\author{Pierangelo Marcati}
\address{Dipartimento di Matematica Pura and Applicata\\
Universit\`a di L'Aquila\\
67100 L'Aquila, Italy}
\email{marcati@univaq.it}

\begin{abstract}
This paper is concerned with the rigorous  analysis of the  zero electron mass limit of the  full Navier-Stokes-Poisson. This system has been introduced in the literature by Anile and Pennisi (see [5]) in order to describe a hydrodynamic model for charge-carrier transport in semiconductor devices. The purpose of this paper is to prove rigorously zero electron mass limit in the framework of general ill prepared initial data. In this situation the velocity field and the electronic fields  develop fast oscillations in time. The main idea we will use in this paper is a combination of formal asymptotic expansion and rigorous uniform estimates on the error terms. Finally we prove the strong convergence of the full Navier Stokes Poisson system towards the incompressible Navier Stokes equations.
\end{abstract}
\maketitle

\section{Introduction}

\subsection{Model}

In this paper we are concerned with the rigorous asymptotic analysis of the following scaled full Navier-Stokes-Poisson type system
\begin{eqnarray}
&&n_t+{\rm div} (nv)=0,\label{s1} \\
&&n[v_t+v\cdot\nabla v] +\nabla p_{m_e}={\rm div}(\frac{\lambda}{m_e} \mathbb{S}(v))+\dfrac{n}{m_e}q \nabla V -\dfrac{nv}{\tau_p},\label{s2}\\
&&\partial_t (nk_BT) +v\cdot\nabla(nk_BT)+\dfrac{5}{3}(nk_BT){\rm div} v +\dfrac{1}{3}{\rm div} q_{m_e} +\dfrac{2}{3}(\lambda\mathbb{S}(v):\nabla v) \nonumber\\
&&\hspace{2cm}=\dfrac{2}{3}n m_e [\dfrac{1}{\tau_p}-\dfrac{1}{2\tau_w}]|v|^2 +\dfrac{m_e n}{2\tau_w}k_B (T_0-T),\label{s3}\\
&&\Delta V=n-\bar n.\label{s4}
\end{eqnarray}
Our goal is to investigate rigorously  the incompressible type limits arising when $m_e$ tends to zero. This kind of incompressible limits are subject to various physical interpretations, the most immediate one, in some analogy with perturbative methods in quantum field theory, is related to the so called {\it zero-mass limit}  in plasma physics (e.g. \cite{GJP99}, \cite{JP00}, \cite{ACJP10}, \cite{AC11}, \cite{CCZ11} ) or otherwise we may consider the quasineutral type limit arising when the Debye length is of the same order of the Mach Number \cite{DM08}. Other similar limits have been investigated by \cite{DM12},  \cite{GM01a}, \cite{GM01b}.

The system  \eqref{s1}-\eqref{s4} has been introduced in the literature by Anile and Pennisi (see \cite{AP92}) in order to  describe a hydrodynamic model for charge-carrier transport in semiconductor devices. Here is the list of notations within \eqref{s1}-\eqref{s4},
\begin{eqnarray*}
\lambda &= n k_B T\tau_\sigma, &\qquad
q_{m_e}=-\dfrac{5}{2m_e}k_B^2nT\tau_q\nabla T +\dfrac{5}{2}(nk_BT)[\dfrac{1}{\tau_p}-\dfrac{1}{\tau_q}]\tau_qv,\\
p_{m_e}&=\dfrac{nk_BT}{m_e} &\qquad
\mathbb{S}(v)=\nabla v+(\nabla v)^T-\dfrac{2}{3}{\rm div} v \mathbb{I},
\end{eqnarray*}
where  $m_{e}$ is the effective electron mass, $k_{B}$ the Boltzmann constant and
$\tau_{p}$, $\tau_{w}$,  $\tau_{q}$ are respectively the momentum, energy,  total energy flow vector relaxation time.
The aim of these models is to incorporate higher-order effects than those included in the standard drift-diffusion equations, in order to be able to describe high-field transport phenomena in semiconductors.
These models includes the fundamental balance laws of particle density, momentum and energy for the
charge carriers and are derived from the moment equations of the Boltzmann transport equation (BTE).

In particular the right-hand sides of \eqref{s1}-\eqref{s4} representing the production of
particles, momentum and energy, due to various interaction mechanisms (carrier-phonon, carrier-carrier and carrier-impurity collisions) are modeled by using  relaxation type nonlinearities.
Moreover in the Anile Pennisi (see \cite{AP92}) derivation they ignore possible anisotropy of the stress tensor and they assume a Fourier-like constitutive law for the heat flux.
The equations of momenta are closed by means of the  principles of extended thermodynamics theory due to Mueller and Ruggeri \cite{MR93}, the hydrodynamic systems then follows by the Maxwellian iteration method.   This choice was motivated by them because of the success of the extended thermodynamics in obtaining the usual system provided by the Grad method for dilute gases. Such an approach allows to generate constitutive laws for the higher momenta of the distribution function, depending on the order of truncation, somehow independently from detailed microscopic assumptions imposed on the form of the distribution function.

The mathematical analysis carried out in this paper gets its inspiration from the papers of \cite{Schochet}, \cite{Schochet94},   \cite{MS01},  \cite{MasmoudiVP},  \cite{ACJP10}, where closed rigorous higher order expansions are carried out in order to manage the different oscillating components contributing to the formal asymptotics.
As far it concerns the Low Mach number limit analysis there is a very large amount of relevant literature that we cannot report here in a reasonable way. Beyond the  classical papers of \cite{KM82}, \cite{Schochet86}, \cite{I87}, \cite{I89}, \cite{U86} we mention  the recent papers of  \cite{A06} where the analysis presents certain similarity with our methods and the review papers of \cite{A08}.
General references on the Navier Stokes equations are for instance the books of  \cite{LPL96, LPL98}, \cite{F04}.

According to the relation among all relaxation constants in Anile and Pennisi's derivation, see \cite{AP92},
\begin{eqnarray*}
\tau_p,\tau_q,\tau_\sigma>0,\quad 2\tau_w>\tau_p,\quad \left[\dfrac{1}{\tau_p}-\dfrac{1}{\tau_q}\right]^2 -\dfrac{4}{5\tau_q}\left[\dfrac{2}{\tau_q}-\dfrac{1}{\tau_w}\right]<0,
\end{eqnarray*}
we can take $\tau_\sigma=\tau_p=\tau_q=\tau_w=1$.

Other physical constants, except the effective mass of electron $m_e$, in system \eqref{s1}-\eqref{s4} will not play any role in our analysis here, we can simply choose $k_B=q=\bar n=1$, and arrive at the following system
\begin{eqnarray*}
&& n_t+{\rm div} (nv)=0,\\
&& n[v_t+v\cdot\nabla v] +\dfrac{1}{m_e}\nabla (n T)=\dfrac{1}{m_e}{\rm div}(nT\mathbb{S}(v)) +\dfrac{n}{m_e} \nabla V-nv,\\
&&\partial_t (nT) +v\cdot\nabla(nT)+\dfrac{5}{3}(nT){\rm div} v -\dfrac{1}{3}\dfrac{5}{2m_e}{\rm div} (nT\nabla T) +\dfrac{2}{3}(nT\mathbb{S}(v):\nabla v,)\\
&&\hspace{2cm}=\dfrac{m_e}{3}n|v|^2 +\dfrac{m_e}{2}n(T_0-T), \\
&&\Delta V=n-1.
\end{eqnarray*}
Moreover, for the convenience of analysis, let $m_e=\varepsilon^2$, our system is changed into
\begin{eqnarray*}
&& n_t+v\cdot \nabla n+ n\nabla\cdot v=0,\\
&& v_t+v\cdot\nabla v +\dfrac{T}{\varepsilon^2 n}\nabla n +\dfrac{1}{\varepsilon^2}\nabla T=\dfrac{1}{\varepsilon^2 n}{\rm div}(nT\mathbb{S}(v)) +\dfrac{1}{\varepsilon^2} \nabla V-v,\\
&& T_t +\dfrac{2}{3}T\nabla\cdot v+v\cdot\nabla T=\dfrac{5}{6\varepsilon^2 n}{\rm div} (nT\nabla T) -\dfrac{2}{3}(T\mathbb{S}(v):\nabla v) +\dfrac{\varepsilon^2}{3}|v|^2 +\dfrac{\varepsilon^2}{2}(T_0-T), \\
&&\Delta V=n-1.
\end{eqnarray*}
This hyperbolic system can be rewritten by using the following notations. Let $U=(n,v,T)^T$ and $F$ be the terms on the right hand side.
$$
U_t+\dsum^d_{j=1}A_j\partial_{x_j}U =F,
$$
where the matrix $A_j$ is defined by
$$
A_j=\left(\begin{array}{ccc} 0 & ne_j & 0 \\ \frac{T}{\varepsilon^2 n}e_j& 0 & \frac{1}{\varepsilon^2}e_j\\ 0 & \frac{2}{3}Te_j & 0\end{array}\right) +v_j I.
$$
It is easy to find that the multiplier to symmetrize the system is
$$
A_0=\left(\begin{array}{ccc}  \frac{T}{\varepsilon^2 n^2}& 0 &0\\ 0 & 1 & 0\\ 0 & 0& \frac{3}{2\varepsilon^2 T}\end{array}\right), \mbox{ such that }
A_0 A_j= \left(\begin{array}{ccc}  \frac{T}{\varepsilon^2 n^2}v_j & \frac{T}{\varepsilon^2 n}e_j& 0 \\[2mm] \frac{T}{\varepsilon^2 n}e_j & v_j & \frac{1}{\varepsilon^2}e_j\\[2mm] 0 & \frac{1}{\varepsilon^2}e_j & \frac{3}{2\varepsilon^2 T}v_j\end{array}\right).
$$

\subsection{Scaled system}

A natural scaling of temperature $T$  should be the one of the same order as electron mass. More precisely, $T$ is of the order $O(\varepsilon^2)$ such that $A_0$ is uniformly positive. Now we use the new variable $\bar T=T/\varepsilon^2$, and the same scaling for $\bar T_0=T_0/\varepsilon^2$, consequently we scale the  electronic potential as $\bar V=V/\varepsilon$. Then the new system is
\begin{eqnarray*}
&& n_t+v\cdot \nabla n+ n\nabla\cdot v=0,\\
&& v_t+v\cdot\nabla v +\dfrac{\bar T}{n}\nabla n +\nabla \bar T=\dfrac{1}{n}{\rm div}(n\bar T\mathbb{S}) +\dfrac{1}{\varepsilon} \nabla \bar V-v,\\
&& \bar T_t +\frac{2}{3}\bar T\nabla\cdot v+v\cdot\nabla \bar T=\dfrac{5}{6 n}{\rm div} (n\bar T\nabla\bar T) -\dfrac{2}{3}(\bar T\mathbb{S}:\nabla v) +\dfrac{1}{3}|v|^2 +\dfrac{\varepsilon^2}{2}(\bar T_0-\bar T), \\
&&\varepsilon\Delta \bar V=n-1.
\end{eqnarray*}

In the remaining discussion, for the convenience of our analysis, we will use the following group of notations, $\sigma^\varepsilon=\dfrac{n-1}{\varepsilon}$, $u^\varepsilon=v$, $T^\varepsilon=\bar T$, $T_0=\bar T_0$, $\psi^\varepsilon=V$, and the system with new notations reads
\begin{align}
&\nonumber \sigma^\varepsilon_t+\dfrac{1}{\varepsilon}{\rm div}u^\varepsilon+ {\rm div}(\sigma^\varepsilon u^\varepsilon)=0\\
&\nonumber u^\varepsilon_t-\dfrac{1}{\varepsilon}\nabla\psi^\varepsilon +u^\varepsilon\cdot\nabla u^\varepsilon = -\dfrac{\varepsilon T^\varepsilon}{1+\varepsilon\sigma^\varepsilon}\nabla \sigma^\varepsilon -\nabla T^\varepsilon\\
&\nonumber\hspace{4cm} +\dfrac{1}{1+\varepsilon\sigma^\varepsilon}{\rm div}((1+\varepsilon\sigma^\varepsilon) T^\varepsilon\mathbb{S}(u^\varepsilon))-u^\varepsilon,\\
&\label{scaledsystem}\tag{$S_{\varepsilon}$} T^\varepsilon_t +\dfrac{2}{3} T^\varepsilon\nabla\cdot u^\varepsilon+u^\varepsilon\cdot\nabla T^\varepsilon = \dfrac{5}{6(1+\varepsilon\sigma^\varepsilon)}{\rm div} ((1+\varepsilon\sigma^\varepsilon) T^\varepsilon\nabla T^\varepsilon) \\
&\nonumber\hspace{4,5cm}-\dfrac{2}{3}( T^\varepsilon\mathbb{S}(u^\varepsilon):\nabla u^\varepsilon) +\dfrac{1}{3}|u^\varepsilon|^2 +\dfrac{\varepsilon^2}{2}(T_0-T^\varepsilon) \\
&\nonumber \Delta  \psi^\varepsilon =\sigma^\varepsilon.
\end{align}
where $\mathbb{S}(u^\varepsilon)=\nabla u^\varepsilon+(\nabla u^\varepsilon)^T-\dfrac{2}{3}{\rm div}u^\varepsilon \mathbb{I}$.

The problem will be studied in a $d$ dimensional periodic domain $\mathbb{T}^d$, with initial data
\begin{eqnarray}
\label{scaledsystemIC}\sigma^\varepsilon|_{t=0}=\sigma_{\rm I}^\varepsilon(x),\quad u^\varepsilon|_{t=0}=u_{\rm I}^\varepsilon(x),\quad T^\varepsilon|_{t=0}=T_{\rm I}^\varepsilon(x).
\end{eqnarray}

The main goal of this paper is to study the limit  as $\varepsilon\rightarrow 0$ in \eqref{scaledsystem}. By a
 simple formal analysis, if we  put $\varepsilon=0$ in \eqref{scaledsystem}, we will get that $n=1$, ${\rm div} v=0$ and our limiting system will be an incompressible Navier-Stokes system with temperature of the form
\begin{eqnarray}
&& {\rm div} v=0,\nonumber\\
&& v_t+ v\cdot\nabla v +\nabla \Pi +\nabla T -{\rm div} (T\mathbb{S}(v)) +v=0,\label{incompressibleNSintro}\\
&& T_t +v\cdot\nabla T =\dfrac{5}{6}{\rm div}(T\nabla T)-\dfrac{2}{3}(T\mathbb{S}(v):\nabla v)+\dfrac{1}{3}|v|^2.\nonumber
\end{eqnarray}
We will give some detailed result on this system in Section \ref{subsecincom}.

The purpose of this paper is to prove rigorously this formal limit in the framework of general ill prepared initial data for the system \eqref{scaledsystem}. In this situation the velocity field and the electronic fields as $\varepsilon \rightarrow 0$ develop fast oscillations in time. In order to describe them one has to introduce
 a linear oscillating limiting system that depends on the solution of \eqref{incompressibleNSintro}.  We  will give a detailed result on the solution of this fast oscillating system in Section \ref{subsecfastosci}.

The main idea we will use in this paper is a combination of formal asymptotic expansion and rigorous uniform estimates on the error terms. The ansatz of our solution is
\begin{eqnarray}
\nonumber\sigma^\varepsilon(x,t) &=& \sigma_{\rm 1f}(x,\frac{t}{\varepsilon}) +\varepsilon \big[\sigma_{\rm 2f}(x,\frac{t}{\varepsilon}) +\sigma_{\rm E}(x,t)\big],\\
\label{expansionintro}u^\varepsilon(x,t) &=& v(x,t)+u_{\rm 1f}(x,\frac{t}{\varepsilon}) +\varepsilon \big[u_{\rm 2f}(x,\frac{t}{\varepsilon}) +u_{\rm E}(x,t)\big],\\
\nonumber T^\varepsilon(x,t) &=& T(x,t)+ \varepsilon \big[T_{\rm 2f}(x,\frac{t}{\varepsilon})+ T_{\rm E}(x,t)\big].
\end{eqnarray}

Here we used notations for subscript ``${\rm 1f}$'' to represent the first order fast oscillation of the solution, ``${\rm 2f}$'' to be the second order fast oscillation of the solution, and ``${\rm E}$'' to be the error terms.

In the above expansion, $(v,T)$ is the solution of limiting incompressible system \eqref{incompressibleNSintro}, $(\sigma_{\rm 1f}, u_{\rm 1f})$ is the solution of first order oscillation system and $(\sigma_{\rm 1f}, u_{\rm 2f}, T_{\rm 2f})$ is the solution of second order oscillation system. The exact form and the solvability of these systems will be listed in the coming discussions.

Since we have the relation $\sigma^\varepsilon=\Delta\psi^\varepsilon$, we have similar representation for electronic potential
\begin{eqnarray*}
\nabla\psi^\varepsilon (x,t) &=& \nabla\psi_{\rm 2f}(x,\dfrac{t}{\varepsilon}) +\varepsilon [\nabla\psi_{\rm 2f}(x,\dfrac{t}{\varepsilon}) +\nabla\psi_{\rm E}(x,t)].
\end{eqnarray*}

For consistence, we also write down initial data in the same form.
\begin{eqnarray*}
\sigma^\varepsilon_{\rm I}(x)&=&\Delta\psi_{\rm I}(x) +\varepsilon\sigma_{\rm I}^E(x),\\
  u^\varepsilon_{\rm I}(x)&=&v_{\rm I}(x)+Qu_{\rm I}(x) +\varepsilon u_{\rm I}^E(x),\\
  T^\varepsilon_{\rm I}(x)&=&T_{\rm I}(x)+\varepsilon T_{\rm I}^E(x),
\end{eqnarray*}
%
%In this paper, we will use the classical notation on the decomposition of $L^2$ functions,
%$$
%u=Pu+Qu.
%$$
where $Q$ is the Leray's projector on the space gradient of vector field $u\in L^{2}(\mathbb{T}^d)$ and is defined as follows,
$$Qu=\nabla \Delta^{-1}{\rm div}u, \qquad Pu=(I-Q)u, \qquad  {\rm div}Pu=0.$$
So, our task is first to find the limiting system, which will be the incompressible system \eqref{incompressibleNSintro} with solution $(v,T,\Pi)$ and a fast oscillation system with solution $(u_{\rm 1f},\sigma_{\rm 1f}=\Delta\psi_{\rm 1f})$. Under suitable assumptions on initial data $v_{\rm I}(x),T_{\rm I}(x)$, ${\rm div}v_{\rm I}=0$, we can prove that \eqref{incompressibleNSintro} is solvable up to time $\tau^*$. Moreover, for any $s>\frac{d}{2}+1$, $\tau\in (0,\tau^*)$, $(v,T)\in L^\infty([0,\tau];H^{s+3})$ and $\Pi\in L^\infty([0,\tau];H^{s+3})$, $\partial_t\Pi\in L^\infty([0,\tau];H^{s+2})$.
We can also prove that the solution of fast oscillation system satisfies $(\sigma_{\rm 1f},u_{\rm 1f})\in L^\infty([0,\tau];H^{s+3})$, $\nabla\psi_{\rm 1f}\in L^\infty([0,\tau];H^{s+3})$.

Now we are ready to state the main result of this paper.

\subsection{Main result}

\begin{thm} Let be $s>\frac{d}{2}+1$. Assume that the initial conditions \eqref{scaledsystemIC} satisfy
\begin{eqnarray*}
\|\Delta\psi_{\rm I}\|_{H^{s+1}}+\|Qu_{\rm I}\|_{H^{s+1}}+\|\varepsilon\sigma_{\rm I}^E\|_{H^s}+\|u_{\rm I}^E\|_{H^s}+\|T_{\rm I}^E\|_{H^s}\leq C.
\end{eqnarray*}
Then, there exists an $\varepsilon_0>0$, such that for all $ \varepsilon\leq \varepsilon_0$, problem \eqref{scaledsystem} \eqref{scaledsystemIC} has a unique classical solution $(\sigma^\varepsilon,u^\varepsilon,T^\varepsilon)\in L^\infty([0,\tau];H^s)$, $u^\varepsilon,T^\varepsilon\in L^\infty([0,\tau];H^s)\cap L^2([0,\tau];H^{s+1})$ and
\begin{eqnarray*}
&\dsup_{0\leq t\leq \tau}\big[\|\varepsilon(\sigma^\varepsilon-\sigma_{\rm 1f})\|_{H^s}+\|u^\varepsilon-v-u_{\rm 1f}\|_{H^s} +\|T^\varepsilon -T\|_{H^s} +\|\nabla\psi^\varepsilon-\nabla\psi_{\rm 1f}\|_{H^s}\big] \leq C\varepsilon.&\\
&\|u^\varepsilon-v-u_{\rm 1f}\|_{L^2([0,\tau];H^{s+1})} +\|T^\varepsilon -T\|_{L^2([0,\tau];H^{s+1})} \leq C\varepsilon.
\end{eqnarray*}
\label{MT}
\end{thm}

Our paper is arranged as follows. In section 2, we will discuss the limiting incompressible system and the leading order oscillating system and give the results on their existence. In section 3, we use asymptotic expansion to find a second order oscillating system and give the derivation of the error system. In section 4, a detailed proof on existence of solution to the error system by using energy method will be given. In the appendix, we give a motivation on how to get our leading order oscillation system.

\section{Limiting system}
\subsection{Limiting incompressible system}\label{subsecincom}

In this section we give some more details on our limiting incompressible Navier-Stokes system, namely
\begin{eqnarray}
&& {\rm div} v=0,\nonumber\\
&& v_t+ v\cdot\nabla v +\nabla \Pi +\nabla T -{\rm div} (T\mathbb{S}(v)) +v=0,\label{incompressibleNS}\\
&& T_t +v\cdot\nabla T =\dfrac{5}{6}{\rm div}(T\nabla T)-\dfrac{2}{3}(T\mathbb{S}(v):\nabla v)+\dfrac{1}{3}|v|^2.\nonumber
\end{eqnarray}
with initial data
\begin{eqnarray}
\label{incompressibleNSIC} v|_{t=0}=v_{\rm I}(x), \quad {\rm div} v_{\rm I}=0,\quad T|_{t=0}=T_{\rm I}(x)\geq T_{\rm L}>0.
\end{eqnarray}
For the system \eqref{incompressibleNS} the following existence results for smooth solutions  holds.
\begin{thm}\label{thmincom}
If $v_{\rm I},T_{\rm I}\in H^{s+3}$ and $T_{\rm I}\geq T_{\rm L}>0$, then there exists $\tau^*>0$ such that \eqref{incompressibleNS} \eqref{incompressibleNSIC} has a unique solution such that for any $\tau<\tau^*$ the following holds
\begin{eqnarray}
\nonumber&&\dsup_{0\leq t\leq \tau}\big[\|v(\cdot,t)\|_{H^{s+3}} +\|T(\cdot,t)\|_{H^{s+3}}\big]
+\|v,T\|_{L^2([0,\tau];H^{s+4})}\\
\label{incompressibleNSest}&&\hspace{5cm}\leq C(\tau,T_{\rm L}) (\|v_{\rm I}\|_{H^{s+3}}+\|T_{\rm I}\|_{H^{s+3}}).\\
\label{estPi} &&\dsup_{0\leq t\leq \tau}\big[\|\Pi(\cdot,t)\|_{H^{s+3}} +\|\partial_t\Pi(\cdot,t)\|_{H^{s+2}}\big] \leq C(\tau,T_{\rm L}).
\end{eqnarray}
\end{thm}

The proof can be  obtained by extending the result by  Kato \cite{K72} for the incompressible Navier Stokes system to the full system and  by using the parabolic theory for the temperature equation. For completeness we mention here the following existence results regarding the coupling of the incompressible Navier Stokes equations and the  temperature balance equation,  \cite{DP08}, \cite{BS11}, \cite{BMR09}, \cite{BM09}.

\subsection{Leading order oscillation system}
\label{subsecfastosci}
By following the same line of arguments as Masmoudi's \cite{MasmoudiVP} and Schochet's \cite{Schochet}  we  describe the main fast oscillating system as follows,
\begin{eqnarray}
\label{oscqlim}
2\partial_t\nabla q &= &Q\left\{-(\nabla q\cdot\nabla)v -(v\cdot\nabla)\nabla q + v\Delta q +{\rm div} [T \mathbb{S}(\nabla q)]\right\}-\nabla q,\\
\label{oscphilim}
2\partial_t\nabla \phi &= &Q\left\{-(\nabla \phi\cdot\nabla)v -(v\cdot\nabla)\nabla \phi + v\Delta \phi+ {\rm div} [T \mathbb{S}(\nabla \phi)]\right\}-\nabla \phi.
\end{eqnarray}
with initial data
\begin{eqnarray}
\label{oscIClim}\nabla q |_{t=0}=Q u_{\rm I}(x),\quad \nabla \phi|_{t=0}=\nabla\psi_{\rm I}(x),
\end{eqnarray}
where $v$ is a divergence free vector field. The details on the construction of such a system can be found in the appendix.

Then,  our leading order fast oscillating system can be obtained from \eqref{oscqlim} and \eqref{oscphilim} by setting
$$
\left(\begin{array}{c} u_{\rm 1f}\\ \nabla\psi_{\rm 1f}\end{array}\right) =e^{-\frac{t}{\varepsilon}L} \left(\begin{array}{c}\nabla q\\ \nabla\phi\end{array}\right), \quad\mbox{ and } \quad \sigma_{\rm 1f}=\Delta\psi_{\rm 1f},
$$
where for any $v,e \in L^2(\T^d;\R^d)$ with ${\rm div}v={\rm div } e=0$, $L$ is defined as follows
\begin{eqnarray*}
& L\left(\begin{array}{c} v\\0\end{array}\right)=L\left(\begin{array}{c} 0\\ e\end{array}\right)=0,\quad&\\
&\mbox{ and }\quad L\left(\begin{array}{c} \nabla q\\ \nabla\phi\end{array}\right)
=\left(\begin{array}{cc} 0& -\I\\ \I & 0\end{array}\right)
\left(\begin{array}{c} \nabla q\\ \nabla \phi\end{array}\right).&
\end{eqnarray*}

Accordingly the system for the fast oscillating vector fields $(\sigma_{\rm 1f},u_{\rm 1f})$ is given by
\begin{align}
\label{rhoosc1f} \partial_t\sigma_{\rm 1f} +\dfrac{1}{2}{\rm div}\big\{ (\nabla\psi_{\rm 1f}\cdot\nabla)v + (v\cdot\nabla)\nabla\psi_{\rm 1f}-&v\sigma_{\rm 1f}
-{\rm div} (T\mathbb{S}(\nabla\psi_{\rm 1f}))\big\}
\\&\nonumber
+\sigma_{\rm 1f} +\dfrac{1}{\varepsilon}{\rm div} u_{\rm 1f} =0,
\end{align}
\begin{align}
\label{uosc1f}\partial_t u_{\rm 1f} +\dfrac{1}{2}Q \big\{ (u_{\rm 1f}\cdot\nabla )v +(v\cdot\nabla) u_{\rm 1f} -v\nabla\cdot u_{\rm 1f} -{\rm div} (T\mathbb{S}(u_{\rm 1f}))\big\} +u_{\rm 1f}-\dfrac{1}{\varepsilon}\nabla\psi_{\rm 1f}=0,
\end{align}
and the equation for $\nabla\psi_{\sigma_{\rm 1f}}$ is
\begin{eqnarray}
\nonumber\partial_t \nabla\psi_{\rm 1f} +\dfrac{1}{2}Q \big\{ (\nabla\psi_{\rm 1f}\cdot\nabla )v +(v\cdot\nabla) \nabla\psi_{\rm 1f} -v\Delta\psi_{\rm 1f} -{\rm div} (T\mathbb{S}(\nabla \psi_{\rm 1f}))\big\} +\nabla\psi_{\rm 1f}+\dfrac{1}{\varepsilon}u_{\rm 1f}=0.
\end{eqnarray}
The initial data of this oscillation system are
\begin{eqnarray}
\label{osc1fIC}
\nabla\psi_{\rm 1f}|_{t=0}=\nabla\psi_{\rm I},\quad u_{\rm 1f}|_{t=0}= Qu_{\rm I},\quad \sigma_{\rm 1f}=\Delta\psi_{\rm I}.
\end{eqnarray}

We  get the existence of the leading order oscillation system immediately from Lemma \ref{lemexistenceqphi} in the appendix.
\begin{lem} \label{lem1f}
The Cauchy problem \eqref{uosc1f}\eqref{rhoosc1f}\eqref{osc1fIC} has a unique solution $(\sigma_{\rm 1f},u_{\rm 1f},\nabla\psi_{\rm 1f})$ which satisfies the following uniform in $\varepsilon$ estimates
\begin{eqnarray}
\label{osc1fest}\dsup_{0\leq t\leq \tau} \big[\|\sigma_{\rm 1f}\|_{H^{s+3}} +\|u_{\rm 1f}\|_{H^{s+3}} +\|\nabla\psi_{\rm 1f}\|_{H^{s+4}}\big] \leq C(\tau).
\end{eqnarray}
\end{lem}

\section{Expansion of the system \eqref{scaledsystem}}

\subsection{Second order fast oscillation system}
By direct calculations from the expansion \eqref{expansionintro}, and by comparing the system with a combination of its incompressible system and the leading order oscillation, we can get that the second order fast oscillations  $(\sigma_{\rm 2f},u_{\rm 2f},T_{\rm 2f},\nabla\psi_{\rm 2f})$ satisfy the following system,
\begin{eqnarray}
\nonumber\partial_s \sigma_{\rm 2f} +{\rm div} u_{\rm 2f}&=&F_{\sigma}^{\rm 2f},\\
\nonumber\partial_s u_{\rm 2f} -\nabla \psi_{\rm 2f} &=& F_{u}^{\rm 2f},\\
\label{osc2f}\partial_s T_{\rm 2f} &=& F_{T}^{\rm 2f},\\
\nonumber\Delta\psi_{\rm 2f}&=&\sigma_{\rm 2f},
\end{eqnarray}
with initial data
\begin{eqnarray}
\label{osc2fIC}\sigma_{\rm 2f}|_{s=0}=u_{\rm 2f}|_{s=0}&=&T_{\rm 2f}|_{s=0}=0,
\end{eqnarray}
where the right hands of the system is given by
\begin{eqnarray*}
F_{\sigma}^{\rm 2f}&=&\dfrac{1}{2}{\rm div} \big\{(\nabla\psi_{\rm 1f}\cdot\nabla)v + (v\cdot\nabla )\nabla\psi_{\rm 1f} -v\sigma_{\rm 1f} -{\rm div} (T\mathbb{S}(\nabla\psi_{\rm 1f}))\big\} \\
&&\hspace{3cm}+\sigma_{\rm 1f}- {\rm div}(\sigma_{\rm 1f}(v+u_{\rm 1f})),\\
F_{u}^{\rm 2f} &=& \dfrac{1}{2}Q \big\{(u_{\rm 1f}\cdot\nabla)v +(v\cdot\nabla) u_{\rm 1f} -v\nabla\cdot u_{\rm 1f} -{\rm div}(T\mathbb{S}(u_{\rm 1f}))\big\} \\
&&\hspace{3cm} -u_{\rm 1f}\nabla (v+ u_{\rm 1f}) -v\cdot\nabla u_{\rm 1f} +{\rm div} (T\mathbb{S}(u_{\rm 1f})),\\
F_T^{\rm 2f}&=& -\dfrac{2}{3}T\nabla u_{\rm 1f} - u_{\rm 1f}\nabla T -\dfrac{2}{3}T\mathbb{S}(v)\nabla u_{\rm 1f} \\
&&\hspace {3cm}-\dfrac{2}{3}T\mathbb{S}(u_{\rm 1f}):\nabla (v+u_{\rm 1f}) -\dfrac{1}{3}|v|^2 +\dfrac{1}{3}|v+u_{\rm 1f}|^2.
\end{eqnarray*}

For the linear system \eqref{osc2f}, one  easily proves the following result
\begin{lem}\label{lem2f}
Problem \eqref{osc2f}\eqref{osc2fIC} has a unique solution %$(\sigma_{\rm 2f},u_{\rm 2f},\nabla\psi_{\rm 2f},T_{\rm 2f})$
which satisfies the following uniform in $\varepsilon$ estimates
\begin{eqnarray}
\label{osc2fest}\dsup_{0\leq s\leq \infty} \big[\|\sigma_{\rm 2f}(\cdot,s)\|_{H^{s+1}} +\|u_{\rm 2f}(\cdot,s)\|_{H^{s+1}} +\|\nabla\psi_{\rm 2f}(\cdot,s)\|_{H^{s+2}}+\|T_{\rm 2f}(\cdot,s)\|_{H^{s+1}}\big] \leq C(\tau).
\end{eqnarray}
\end{lem}

\subsection{$O(\varepsilon)$ Approximation system}

Summing up the limiting incompressible system \eqref{incompressibleNS}, the first order oscillation system \eqref{rhoosc1f}\eqref{uosc1f} and the second order oscillation system \eqref{osc2f}, we have that the $O(\varepsilon)$ correction of our scaled system \eqref{scaledsystem} is
\begin{align}
\nonumber &\partial_t (\sigma_{\rm 1f}+\varepsilon \sigma_{\rm 2f})+\dfrac{1}{\varepsilon}{\rm div} (u_{\rm 1f}+\varepsilon u_{\rm 2f}) +{\rm div} (\sigma_{\rm 1f} (v+u_{\rm 1f}))=0,\\
\nonumber & \partial_t (v+u_{\rm 1f}+\varepsilon u_{\rm 2f}) -\dfrac{1}{\varepsilon}\nabla(\psi_{\rm 1f}+\varepsilon \psi_{\rm 2f})+(v+u_{\rm 1f})\\
\label{appsystem} \tag{$S_{1f}$}&\hspace{2cm} +(v+u_{\rm 1f})\cdot \nabla (v+u_{\rm 1f}) -{\rm div} (T \mathbb{S}(v+u_{\rm 1f})) +\nabla\Pi +\nabla T=0,\\
\nonumber & \partial_t (T+\varepsilon T_{\rm 2f}) +\dfrac{2}{3}T\nabla\cdot (v+u_{\rm 1f}) +(v+u_{\rm 1f})\cdot \nabla T \\
\nonumber &\hspace{2cm}+\dfrac{2}{3}T\mathbb{S}(v+u_{\rm 1f}) :\nabla (v+u_{\rm 1f}) -\dfrac{5}{6}{\rm div} (T\nabla T) -\dfrac{1}{3} |v+u_{\rm 1f}|^2=0.
\end{align}
We will use the following notations for  the approximation of the solution.
\begin{eqnarray}
\nonumber\sigma_{\rm app}(x,t)&=& \sigma_{\rm 1f}\left(x,\frac{t}{\varepsilon}\right) +\varepsilon \sigma_{\rm 2f}\left(x,\frac{t}{\varepsilon}\right)\\
\label{appsol}u_{\rm app}(x,t)&=& v(x,t)+u_{\rm 1f}\left(x,\frac{t}{\varepsilon}\right)+\varepsilon u_{\rm 2f}\left(x,\frac{t}{\varepsilon}\right)\\
\nonumber T_{\rm app}(x,t)&=& T(x,t)+\varepsilon T_{\rm 2f}\left(x,\frac{t}{\varepsilon}\right),
\end{eqnarray}

The initial data of our approximation problem is
\begin{eqnarray}
\label{appIC} \sigma_{\rm app}|_{t=0}=\Delta\psi_{\rm I}(x),\quad u_{\rm app}|_{t=0}=v_{\rm I}(x)+Qu_{\rm I}(x), \quad T_{\rm app}|_{t=0}=T_{\rm I}(x).
\end{eqnarray}

Moreover, by the results of Theorem \ref{thmincom}, Lemma \ref{lem1f} and Lemma \ref{lem2f}, we have the following uniform estimates for this approximation solution
\begin{eqnarray}
\label{appsolest}
&\dsup_{0\leq t\leq \tau} \big[\|\sigma_{\rm app}\|_{H^{s+1}} +\|u_{\rm app}\|_{H^{s+1}} +\|\nabla\psi_{\rm app}\|_{H^{s+2}}+\|T_{\rm app}\|_{H^{s+1}}\big] \leq C(\tau).&\\
\nonumber &T_{\rm app}\geq T_L/2>0.&
\end{eqnarray}

\subsection{System for error terms}

Now, by comparing $O(\varepsilon)$ correction system \eqref{appsystem} with our scaled system \eqref{scaledsystem}, we get  the following system for the error terms $(\sigma_{\rm E},u_{\rm E},T_{\rm E},\psi_{\rm E})$,
\begin{align}
\nonumber &\partial_t(\sigma_{\rm E}+\Delta\Pi) +u_{\rm app}\cdot\nabla (\sigma_{\rm E}+\Delta\Pi) +\dfrac{1+\varepsilon\sigma^\varepsilon}{\varepsilon}{\rm div} u_{\rm E}= G_{\sigma}^{\rm E} + F_{\sigma}^{\rm E},\\
\label{errorsystemv2}\tag{$S_{E}$}&\partial_t u_{\rm E} +u_{\rm app}\cdot\nabla u_{\rm E} +\dfrac{\varepsilon T_{\rm app}}{1+\varepsilon\sigma^\varepsilon}\nabla (\sigma_{\rm E}+\Delta\Pi) +\nabla T_{\rm E}+u_{\rm E} -\dfrac{1}{\varepsilon}\nabla (\psi_{\rm E} +\Pi)\\
\nonumber &\hspace{5.7cm} -{\rm div}(T_{\rm app}\mathbb{S}(u_{\rm E}))= G_{u}^{\rm E} + F_{u}^{\rm E},
\\
\nonumber &\partial_t T_{\rm E}+ \dfrac{2}{3}T_{\rm app}\nabla\cdot u_{\rm E} +u_{\rm app}\cdot\nabla T_{\rm E} -\dfrac{5}{6}{\rm div}(T_{\rm app}\nabla T_{\rm E}) =G_T^{\rm E} + F_T^{\rm E},
\end{align}
with initial data
\begin{eqnarray}
\label{errorsystemIC}\sigma_{\rm E}|_{t=0}=\sigma_{\rm I}^E,\quad u_{\rm E}|_{t=0}=u_{\rm I}^E,\quad T_{\rm E}|_{t=0}=T_{\rm I}^E.
\end{eqnarray}
The right hand sides of the error system \eqref{errorsystemv2} are given by
\begin{eqnarray*}
G_{\sigma}^{\rm E}&=& -\nabla(\sigma_{\rm app}+\varepsilon\sigma_{\rm E})\cdot u_{\rm E} -\sigma_{\rm E}{\rm div} u_{\rm app},\\
G_{u}^{\rm E}&=& -u_{\rm E}\cdot\nabla (u_{\rm app}+\varepsilon u_{\rm E}) + {\rm div} (T_{\rm E}\mathbb{S}(u_{\rm app}+\varepsilon u_{\rm E})) +\dfrac{\varepsilon T_{\rm app}}{1+\varepsilon\sigma^\varepsilon}\nabla (\sigma_{\rm E}+\Delta\Pi) \\
&&\hspace{0.5cm}-\dfrac{T^\varepsilon}{1+\varepsilon\sigma^\varepsilon}\nabla\sigma^\varepsilon +\dfrac{\nabla\sigma^\varepsilon}{1+\varepsilon\sigma^\varepsilon} T^\varepsilon \mathbb{S}(u^\varepsilon), \\
G_T^{\rm E} &=& -\dfrac{2}{3} T_{\rm E}\nabla\cdot (u_{\rm app}+\varepsilon u_{\rm E}) - u_{\rm E}\cdot\nabla (T_{\rm app}+\varepsilon T_{\rm E}) +\dfrac{5}{6}{\rm div}(T_{\rm E}\nabla (T_{\rm app}+\varepsilon T_{\rm E})) \\
 &&\hspace{0.5cm} -\dfrac{2}{3}\varepsilon^2 T_{\rm E} \mathbb{S}(u_{\rm E}):\nabla u_{\rm E} -\dfrac{2}{3}\varepsilon T_{\rm E} \mathbb{S}(u_{\rm E}): \nabla u_{\rm app}\\
 &&\hspace{0.5cm}   -\dfrac{2}{3}\varepsilon T_{\rm E} \mathbb{S}(u_{\rm app})  :\nabla u_{\rm E} -\dfrac{2}{3}\varepsilon T_{\rm app} \mathbb{S}(u_{\rm E}):\nabla u_{\rm E}\\
&& \hspace{0.5cm} -\dfrac{2}{3} T_{\rm app} \mathbb{S}(u_{\rm E}) : u_{\rm app} -\dfrac{2}{3} T_{\rm app} \mathbb{S}(u_{\rm app}):\nabla u_{\rm E} -\dfrac{2}{3}T_{\rm E}\mathbb{S}(u_{\rm app}) : \nabla u_{\rm app}\\
 &&\hspace{0.5cm}+\dfrac{\varepsilon}{3}|u_{\rm E}|^2 -\dfrac{\varepsilon^2}{2}T_{\rm E}  +\dfrac{5\nabla\sigma^\varepsilon}{6(1+\varepsilon\sigma^\varepsilon)}T^\varepsilon\nabla T^\varepsilon+\dfrac{2}{3}u_{\rm app} \cdot u_{\rm E},
\end{eqnarray*}
and
\begin{eqnarray*}
F_{\sigma}^{\rm E}&=& -{\rm div} (\sigma_{\rm 2f}u_{\rm app}) -{\rm div} (\sigma_{\rm 1f}u_{\rm 2f}) +\partial_t\Delta\Pi +u_{\rm app}\nabla\Delta\Pi. \\
F_u^{\rm E}&=& -u_{\rm 2f}-\nabla T_{\rm 2f} - u_{\rm 2f}\cdot\nabla u_{\rm app} -(v+u_{1f})\cdot\nabla u_{\rm 2f} +{\rm div}(T_{\rm 2f}\mathbb{S}(u_{\rm app})) +{\rm div}(T\mathbb{S}(u_{\rm 2f})).\\
F_T^{\rm E}&=& -\dfrac{2}{3}  T\nabla\cdot u_{\rm 2f} -\dfrac{2}{3}T_{\rm 2f}\nabla\cdot u_{\rm app} -u_{\rm 2f}\cdot\nabla T_{\rm app} -(v+u_{\rm 1f})\cdot\nabla T_{\rm 2f}\\
&&\hspace{0.5cm} -\dfrac{2}{3}T_{\rm 2f}\mathbb{S}(u_{\rm app}):\nabla u_{\rm app}  -\dfrac{2}{3}T\mathbb{S}(u_{\rm 2f}):\nabla u_{\rm app}  -\dfrac{2}{3}T\mathbb{S}(v+u_{\rm 1f}):\nabla u_{\rm 2f} \\
&&\hspace{0.5cm}+\dfrac{5}{6}{\rm div}(T_{\rm app}\nabla T_{\rm 2f}) +\dfrac{5}{6}{\rm div} (T_{\rm 2f}\nabla T) +\dfrac{2}{3}u_{\rm app}\cdot u_{\rm 2f} +(v+u_{\rm 1f})u_{\rm 2f} +\dfrac{\varepsilon}{2}(T_0-T_{\rm app}).
\end{eqnarray*}
Now, the last step is to prove the following existence result for the error terms.
\begin{prop}\label{properror}
Let $s>\frac{n}{2}+1$. If $\sigma_{\rm I}^E,u_{\rm I}^E,T_{\rm I}^E\in H^s$, then \eqref{errorsystemv2}\eqref{errorsystemIC} has a unique solution such that
\begin{eqnarray}
\label{estError}\dsup_{0\leq t\leq \tau}\|\varepsilon\sigma_{\rm E},u_{\rm E},T_{\rm E},\nabla\psi_{\rm E}\|_{H^s}+ \|u_{\rm E},T_{\rm E}\|_{L^2([0,\tau];H^{s+1})}\leq C.
\end{eqnarray}
\end{prop}
The proof of this proposition will be done in  the next section.

\section{Existence of solution to the system \eqref{errorsystemv2} of the error terms}

In this section, we focus on solving the system of error terms by using energy method.

The following Moser's type inequality will be used several times in the coming estimates. For completeness, we list it here as a proposition.
\begin{prop}\label{propMoser}
The following facts hold
\begin{enumerate}
\item Let $s\ge 0$, $f$, $g\in H^s(\T^d)\cap L^\infty(\T^d)$,
and $\alpha$ a multi-index with $|\alpha|\le s$. Then, for some
constant $c_s>0$,
\begin{equation}
  \|D^\alpha(fg)\|_0 \le c_s(\|f\|_\infty \|D^s g\|_0 +\|g\|_\infty
    \|D^s f\|_0).  \label{moser1}
\end{equation}
\item Let $s\ge 1$, $f\in H^s(\T^d)$ with $Df\in L^\infty(\T^d)$,
$g\in H^{s-1}(\T^d)\cap L^\infty(\T^d)$, and $|\alpha|\le s$. Then,
for some constant $c_s>0$,
\begin{equation}
  \|[D^\alpha,f]g\|_0\leq c_s (\|Df\|_\infty
    \|D^{s-1}g\|_0 +\|g\|_\infty \|D^sf\|_0).  \label{moser2}
\end{equation}
where $[D^\alpha,f]g=D^\alpha(fg)-fD^\alpha g$.
\end{enumerate}
\end{prop}

We rewrite the system \eqref{errorsystemv2} into a symmetrizable hyperbolic system. We introduce some notations to shorten our formula.
\begin{eqnarray*}
U_{\rm E}=\left(\begin{array}{c}\sigma_{\rm E}+\Delta\Pi\\ u_{\rm E} \\ T_{E}\end{array}\right), \quad
A_j^E=\left(\begin{array}{ccc}0&\frac{1+\varepsilon\sigma^\varepsilon}{\varepsilon}e_j & 0\\ \frac{\varepsilon}{1+\varepsilon\sigma^\varepsilon} T_{\rm app}e_j & 0 &e_j\\ 0& \frac{2}{3}T_{\rm app}e_j & 0\end{array}\right) + (u_{\rm app})_j I.
\end{eqnarray*}
Then we can choose the symmetrizer of the system to be
\begin{eqnarray*}
A_0^E=\left(\begin{array}{ccc}\frac{\varepsilon^2 T_{\rm app}}{(1+\varepsilon\sigma^\varepsilon)^2}&0&0\\[2mm] 0&I&0\\[2mm] 0&0& \frac{3}{2T_{\rm app}}\end{array}\right), \quad \mbox{ and }
\mathcal{A}_j^E=A_0^EA_j^E=\left(\begin{array}{ccc} \frac{\varepsilon^2 T_{\rm app}}{(1+\varepsilon\sigma^\varepsilon)^2}(u_{\rm app})_j&\frac{\varepsilon T_{\rm app}}{1+\varepsilon\sigma^\varepsilon}e_j & 0\\[2mm] \frac{\varepsilon T_{\rm app}}{1+\varepsilon\sigma^\varepsilon}e_j & (u_{\rm app})_j &e_j\\[2mm] 0& e_j & \frac{3}{2 T_{\rm app}}(u_{\rm app})_j\end{array}\right).
\end{eqnarray*}

Therefore, the error system \eqref{errorsystemv2} is written down into a symmetric hyperbolic system with viscosity
\begin{eqnarray}
\label{hyperbolic}\partial_t U_{\rm E} +A_j^E \partial_{x_j} U_{\rm E} - \mathcal{D}U_{\rm E}=- \dfrac{1}{\varepsilon} J +H^E,
\end{eqnarray}
where
$$
\mathcal{D}U_{\rm E}=\left(\begin{array}{c} 0\\ -u_{\rm E}+{\rm div}(T_{\rm app}\mathbb{S}(u_{\rm E})) \\ \frac{5}{6}{\rm div}(T_{\rm app}\nabla T_{\rm E}) \end{array}\right),\quad J=\left(\begin{array}{c} 0\\ -\nabla (\psi_{\rm E}+\Pi) \\ 0 \end{array}\right),
$$
and
$$
H^E=H_G^E+H_F^E=\left(\begin{array}{c} G_\sigma^E\\G_u^E\\G_T^E\end{array}\right)+\left(\begin{array}{c} F_\sigma^E\\F_u^E\\F_T^E\end{array}\right).
$$
The symmetric version of our system \eqref{errorsystemv2} is
\begin{eqnarray}
\label{hyperbolicsym}A_0^E\partial_t U_{\rm E} + \mathcal{A}_j^E \partial_{x_j} U_{\rm E} - A_0^E\mathcal{D}U_{\rm E}=- \dfrac{1}{\varepsilon} A_0^EJ +A_0^EH^E,
\end{eqnarray}

In the following we will use the operator
$$
\mathcal{L}U=A_0^E\partial_t U + \mathcal{A}_j^E \partial_{x_j} U - A_0^E\mathcal{D}U,
$$

and, for  any multi-index $\alpha$, we will use the notations
$$
\sigma_\alpha=D^\alpha\sigma_{\rm E},\quad \psi_{\alpha}=D^\alpha (\psi_{\rm E}+\Pi), \quad u_\alpha=D^\alpha u_{\rm E},\quad T_\alpha=D^\alpha T_{\rm E} \quad \mbox{ and } U_\alpha=\left(\begin{array}{c}\sigma_\alpha\\ u_\alpha \\ T_\alpha\end{array}\right),
$$
with the following relation $\sigma_\alpha=\Delta \psi_\alpha$.

Now for any multi-index $\alpha$,  $0\leq |\alpha|\leq s$,  $s>\frac{d}{2}+1$, by applying the differential operator $D^\alpha$ on the hyperbolic system \eqref{hyperbolic} and by symmetrizing it, we have
\begin{eqnarray}
\label{hyperbolicalpha}
&\mathcal{L}U_\alpha=A_0^E\partial_t U_\alpha + \mathcal{A}_j^E \partial_{x_j} U_\alpha - A_0^E\mathcal{D}U_\alpha=- \dfrac{1}{\varepsilon} A_0^ED^\alpha J+\mathcal{R}_\alpha&\\
\label{hyperalphaR}&\mathcal{R}_\alpha=A_0^ED^\alpha H^E + A_0^E[A_j^E, D^\alpha]\partial_{x_j} U_{\rm E} -A_0^E[\mathcal{D},D^\alpha] U_{\rm E},
\end{eqnarray}
where we used the commutator notation $[\cdot,\cdot]$.

By applying the energy method for the system \eqref{hyperbolicalpha}, we obtain
\begin{eqnarray}
\label{energyest}\langle \mathcal{L}U_\alpha, U_\alpha\rangle =-\dfrac{1}{\varepsilon}\dint \nabla\psi_\alpha\cdot u_\alpha +\langle \mathcal{R}_\alpha, U_\alpha\rangle.
\end{eqnarray}
where $\langle f,g\rangle=\dint f g$ and also
$$
\|U_\alpha\|_{\rm e}^2=\dint (|\varepsilon\sigma_\alpha|^2+|u_\alpha|^2+|T_\alpha|^2),\quad \mbox{ and } \|U_{\rm E}\|_{s, {\rm e}} =\dsum^s_{|\alpha|=0} \|U_\alpha\|_{\rm e}.
$$

Now we divide the estimates into three steps. In the remaining estimates, we will denote $C$ to be constants depending only on the following quantities,
\begin{eqnarray*}
&\dsup_{0\leq t\leq\tau}\|U_{\rm E}\|_{s, {\rm e}}, \|\sigma_{\rm 1f},\sigma_{\rm 2f}\|_{L^\infty([0,\tau];H^{s+1})},  \|v, u_{\rm 1f}, u_{\rm 2f}\|_{L^\infty([0,\tau];H^{s+1})}, \\
& \| T,T_{\rm 2f}\|_{L^\infty([0,\tau];H^{s+1})}, \|\Pi\|_{L^\infty([0,\tau];H^{s+3})}, \|\partial_t\Pi\|_{L^\infty([0,\tau];H^{s+2})}, T_{\rm L}.
\end{eqnarray*}

{\bf Step. 1}. The left hand side for linear operator $\mathcal{L}$ in \eqref{energyest} yields to
\begin{eqnarray*}
\langle \mathcal{L}U_\alpha, U_\alpha\rangle &= &\dfrac{1}{2}\dfrac{d}{dt}\dint (\frac{\varepsilon^2 T_{\rm app}}{(1+\varepsilon\sigma^\varepsilon)^2}\sigma_\alpha^2 +|u_\alpha|^2 +\dfrac{3}{2T_{\rm app}}T_\alpha^2) \\
&&\quad +\dint |u_\alpha|^2 +\dint T_{\rm app} \mathbb{S}(u_\alpha):\nabla u_\alpha +\dfrac{5}{4}\dint |\nabla T_\alpha|^2 \\
&&\quad +\dfrac{5}{4}\dint \dfrac{1}{T_{\rm app}}T_\alpha\nabla T_{\rm app}\nabla T_\alpha-\dint (\partial_t A_0^E U_\alpha)\cdot U_\alpha -\dfrac{1}{2}\dint (\partial_{x_j} \mathcal{A}_j^E U_\alpha)\cdot U_\alpha,
\end{eqnarray*}
where the last three terms can be estimated by
\begin{eqnarray*}
C \|U_\alpha\|_{\rm e}^2 +\frac{1}{4}\dint |\nabla T_\alpha|^2.
\end{eqnarray*}

Hence we get
\begin{eqnarray}
\nonumber&&\langle \mathcal{L}U_\alpha, U_\alpha\rangle \geq \dfrac{1}{2}\dfrac{d}{dt}\dint (\frac{\varepsilon^2 T_{\rm app}}{(1+\varepsilon\sigma^\varepsilon)^2}\sigma_\alpha^2 +|u_\alpha|^2 +\dfrac{3}{2T_{\rm app}}T_\alpha^2) \\
\label{calL}&&\hspace{2cm}+\dint |u_\alpha|^2 +\frac{T_l}{2}\dint |\nabla u_\alpha|^2 +\dint |\nabla T_\alpha|^2 - C \|U_\alpha\|_{\rm e}^2.
\end{eqnarray}

{\bf Step. 2.}
The singular term on the right hand side of \eqref{energyest} can be handled by using the first equation of \eqref{hyperbolicalpha}, i.e.
$$
\partial_t \Delta\psi_\alpha +D^\alpha(u_{\rm app}\cdot\nabla \Delta\psi_0) +\dfrac{1}{\varepsilon} {\rm div} u_\alpha +\sigma^\varepsilon {\rm div}u_\alpha = [\frac{1+\varepsilon\sigma^\varepsilon}{\varepsilon},D^\alpha]{\rm div}u_{\rm E}+ D^\alpha (G_\sigma^E+F_\sigma^E).
$$
Now the singular term in the energy estimate \eqref{energyest} is
\begin{eqnarray}
\nonumber&&-\dfrac{1}{\varepsilon}\dint \nabla\psi_\alpha\cdot u_\alpha =\dfrac{1}{\varepsilon}\dint \psi_\alpha {\rm div}u_\alpha\\
\nonumber&=& -\dint \psi_\alpha \partial_t\Delta\psi_\alpha -\dint \psi_\alpha D^\alpha( u_{\rm app}\cdot\nabla \Delta\psi_0) -\dint\psi_\alpha\sigma^\varepsilon {\rm div}u_\alpha\\
\label{singularterm} &&\hspace{2cm}+\dint\psi_\alpha [\frac{1+\varepsilon\sigma^\varepsilon}{\varepsilon},D^\alpha]{\rm div}u_{\rm E}+\dint \psi_\alpha D^\alpha\big(G_\sigma^E+F_\sigma^E\big)\\
\nonumber&=& -\dfrac{1}{2}\dfrac{d}{dt}\dint |\nabla \psi_\alpha|^2  -\dint \psi_\alpha u_{\rm app}\cdot\nabla \Delta\psi_\alpha -\dint \psi_\alpha [D^\alpha, u_{\rm app}] \nabla\Delta\psi_0 -\dint\psi_\alpha\sigma^\varepsilon {\rm div}u_\alpha\\
\nonumber && \hspace{2cm} +\dint\psi_\alpha [\frac{1+\varepsilon\sigma^\varepsilon}{\varepsilon},D^\alpha]{\rm div}u_{\rm E}+\dint \psi_\alpha D^\alpha\big(G_\sigma^E+F_\sigma^E\big),
\end{eqnarray}
where the fist term on the right hand side contributes in the energy.

The last four terms  in \eqref{singularterm} are easily estimated by the energy as we will do later, while the second term is different. Since we have only $\varepsilon\sigma_\alpha =\varepsilon\Delta\psi_\alpha$ in the energy, there is no hope to control it just directly by $\Delta\psi_\alpha$ and $\nabla\psi_\alpha$. Our idea here is to move one of the derivatives onto $u_{\rm app}$ which is known and $W^{2,\infty}$ controllable. The reason we can succeed in doing this is that this term is quadratic in $\psi_\alpha$. More precisely,
\begin{eqnarray*}
&&-\dint \psi_\alpha u_{\rm app} \cdot \nabla\Delta \psi_\alpha=\dint \nabla\psi_\alpha \cdot u_{\rm app}\Delta \psi_\alpha +\dint \psi_\alpha {\rm div} u_{\rm app} \Delta\psi_\alpha\\
&=&\dsum^n_{i,j=1}\dint \partial_i\psi_\alpha u_{\rm app}^i \partial_{jj}\psi_\alpha +\dint \psi_\alpha {\rm div} u_{\rm app} \Delta\psi_\alpha\\
&=& \dsum^n_{i,j=1}\dint \partial_{ij}\psi_\alpha u_{\rm app}^i \partial_{j}\psi_\alpha +\dsum^n_{i,j=1}\dint \partial_i\psi_\alpha \partial_j u_{\rm app}^i \partial_{j}\psi_\alpha -\dint \nabla(\psi_\alpha {\rm div} u_{\rm app}) \nabla\psi_\alpha\\
&=& \dfrac{1}{2}\dsum^n_{i,j=1}\dint u_{\rm app}^i \partial_i(\partial_{j}\psi_\alpha)^2 +\dint (\nabla\psi_\alpha \nabla u_{\rm app})\cdot \nabla\psi_\alpha \\
&&-\dint \nabla\psi_\alpha {\rm div} u_{\rm app} \nabla\psi_\alpha -\dint \psi_\alpha \nabla{\rm div} u_{\rm app} \nabla\psi_\alpha\\
&=& -\dfrac{1}{2}\dint {\rm div} u_{\rm app} |\nabla\psi_\alpha|^2  +\dint (\nabla\psi_\alpha \nabla u_{\rm app})\cdot \nabla\psi_\alpha \\
&&-\dint \nabla\psi_\alpha {\rm div} u_{\rm app} \nabla\psi_\alpha -\dint \psi_\alpha \nabla{\rm div} u_{\rm app} \nabla\psi_\alpha.
\end{eqnarray*}
Thus this term can be estimated by
\begin{eqnarray*}
-\dint \psi_\alpha u_{\rm app} \cdot \nabla\Delta \psi_\alpha &\leq & \|u_{\rm app}\|_{W^{2,\infty}} \dint (|\nabla \psi_\alpha|^2+|\psi_\alpha|^2).
\end{eqnarray*}

The estimate for all the  terms, but the last one,  on the right hand side of \eqref{singularterm} can be done by integral by parts and Moser's type inequalities Proposition \ref{propMoser}, i.e.
\begin{eqnarray*}
\Big|-\dint \psi_\alpha [D^\alpha, u_{\rm app}] \nabla\Delta\psi_0 \Big|\leq \|u_{\rm app}\|_{W^{2,\infty}} \|\psi_0\|_{H^s}^2,
\end{eqnarray*}
\begin{eqnarray*}
\Big|\dint\psi_\alpha\sigma^\varepsilon {\rm div}u_\alpha\Big| \leq C(\|\varepsilon\sigma_0\|_{H^s}^2+\|\varphi_\alpha\|_{L^2}^2) +\dfrac{1}{4}\dint |\nabla u_\alpha|^2,
\end{eqnarray*}
and
\begin{eqnarray*}
\left|\dint\psi_\alpha [\frac{1+\varepsilon\sigma^\varepsilon}{\varepsilon},D^\alpha]{\rm div}u_{\rm E}\right|\leq C (\|\phi_\alpha\|_{L^2}^2+\|\varepsilon\sigma_0\|_{H^s}^2+\|u_{\rm E}\|_{H^s}^2).
\end{eqnarray*}

The last term in \eqref{singularterm} can be estimated by using integral by parts, H\"older's inequality and Moser's type inequalities, Proposition \ref{propMoser},
\begin{eqnarray*}
\Big|\dint \psi_\alpha D^\alpha\big(G_\sigma^E+F_\sigma^E\big)\Big|&\leq  & \|\nabla\psi_\alpha\|_{L^2}(\|D^{\alpha-1}G_\sigma^E\|_{L^2}+\|D^{\alpha-1}F_\sigma^E\|_{L^2})\\
&\leq & C\big(\|\nabla\psi_0\|_{H^s}^2+\|\varepsilon\sigma_0\|_{H^s}^2 +\|u_E\|_{H^s}^2\big).
\end{eqnarray*}
%where $C$ depends on $\|\varepsilon(\sigma_{\rm E}+\Delta\Pi\|_{L^\infty}$, $\|u_{\rm E}\|_{L^\infty}$, $\|\sigma_{\rm app}\|_{H^s}$, $\|u_{\rm %app}\|_{H^{s+1}}$, $\|\sigma_{\rm 1f},\sigma_{\rm 2f}\|_{H^{s-1}}$, $\|u_{\rm 1f}\|_{H^{s-1}}$, $\|\partial_t\Pi\|_{H^{s+1}}$ and $\|\Pi\|_{H^{s+2}}$.

%We need to find out here which kind of regularity we need for this term, the same for the estimates in step 3.

Then our estimates on the singular term in \eqref{energyest} is
\begin{eqnarray}
-\dfrac{1}{\varepsilon}\dint \nabla\psi_\alpha\cdot u_\alpha \leq -\dfrac{1}{2}\dint {\rm div} u_{\rm app} |\nabla\psi_\alpha|^2 +C (\|U_{\rm E}\|_{s,{\rm e}}^2+\|\nabla \psi_0\|_{H^s}^2) +\dfrac{1}{4}\dint |\nabla u_\alpha|^2.
\label{singulartermend}
\end{eqnarray}

{\bf Step. 3.} Finally, we estimate the remainder terms $\langle \mathcal{R}_\alpha, U_\alpha\rangle $ in \eqref{energyest}.

The commutator terms are controlled by using Moser type inequality, Proposition \ref{propMoser},
\begin{eqnarray*}
\left|\dint A_0^E[A_j^E,D^\alpha] \partial_{x_j} U_{\rm E}U_\alpha -\dint A_0^E[\mathcal{D},D^\alpha] U_{\rm E} U_\alpha\right|\leq C (\|u_{\rm E}\|_{H^{s+1}}+\|T_{\rm E}\|_{H^{s+1}}) \|U_\alpha\|_{\rm e},
\end{eqnarray*}

The last error term is
\begin{eqnarray*}
&&\dint A_0^ED^\alpha H^E \cdot U_\alpha =\dint A_0^ED^\alpha H^E_G \cdot U_\alpha+ \dint A_0^ED^\alpha H^E_F \cdot U_\alpha.
\end{eqnarray*}

Since we have the viscosity and the heat diffusion terms in our system, we can move one derivative from the nonlinear term to $u_\alpha$ and $T_\alpha$. While for density $\sigma_\alpha$, we couldn't do this.

Now we know that all $H_F^E$ is given by our approximation function which are known. We can simply estimate the second term above by
\begin{eqnarray*}
&&\Big|\dint A_0^ED^\alpha H^E_F \cdot U_\alpha\Big| \\
&\leq & C\big(\|D^\alpha F_\sigma^E\|_{L^2} \|\varepsilon\sigma_0\|_{H^s} +\|D^{\alpha-1} F_u^E\|_{L^2} \|\nabla u_\alpha\|_{L^2} +\|D^{\alpha-1} F_T^E\|_{L^2}\|\nabla T_0\|_{H^s}\big)\\
&\leq & C+C\|\varepsilon\sigma_0\|_{H^s}^2 +\eta\|\nabla u_\alpha\|_{L^2}^2 +\eta\|\nabla T_0\|_{H^s}^2, \qquad \mbox{ for small }\eta.
\end{eqnarray*}

The estimates on $H_G^E$, $G_u^E$ and $G_T^E$ are relatively easier by using Moser's type inequalities, Proposition \ref{propMoser}. Just notice that whenever we have a $D^{\alpha+1}$ on $u_{\rm E}$ or $T_{\rm E}$, there is an $\varepsilon$ in front, which give us the possibility to control them by using the viscosity and the heat diffusion terms. The result is
\begin{eqnarray*}
&&\Big|\dint D^\alpha G_u^E u_\alpha +\dint D^\alpha G_T^E \dfrac{3}{2T_{\rm app}} T_\alpha\Big|\\
&\leq & \|D^{\alpha-1} G_u^E\|_{L^2} \| \nabla u_\alpha\|_{L^2} + \| D^\alpha G_T^E \|_{L^2} \|\dfrac{3}{2T_{\rm app}} T_\alpha\|_{L^2}\\
&\leq & C + C\|U_E\|_{s,{\rm e}}^2 +\|\varepsilon\nabla u_\alpha\|_{L^2}^2 +\|\varepsilon\nabla T_\alpha\|_{L^2}^2.
\end{eqnarray*}

But we have to be more careful with the term $D^\alpha G_\sigma^E$, since we could not move the derivative to $\varepsilon\sigma_\alpha$ since there is $\nabla \sigma_{\rm E}$ in $G_\sigma^E$ which will give us $D^{\alpha+1}\sigma_{\rm E}$ and we will not close our estimates with this. Our method here is similar to the one we used to deal with singular term \eqref{singularterm}. More precisely, we have
\begin{eqnarray*}
&&\Big|\dint D^\alpha G_\sigma^{\rm E} \dfrac{\varepsilon^2T_{\rm app}}{(1+\varepsilon\sigma^\varepsilon)^2} \sigma_\alpha\Big| \\
&\leq &\Big|\dint D^\alpha (\varepsilon\nabla \sigma_{\rm E}\cdot u_{\rm E})\dfrac{\varepsilon^2T_{\rm app}}{(1+\varepsilon\sigma^\varepsilon)^2} \sigma_\alpha\Big|+\Big|\dint D^\alpha (\nabla \sigma_{\rm app}\cdot u_{\rm E} +\sigma_{\rm E}{\rm div}u_{\rm app})\dfrac{\varepsilon^2T_{\rm app}}{(1+\varepsilon\sigma^\varepsilon)^2} \sigma_\alpha\Big|\\
&\leq & \Big|\dint \varepsilon\nabla \sigma_\alpha\cdot u_{\rm E}\dfrac{\varepsilon^2T_{\rm app}}{(1+\varepsilon\sigma^\varepsilon)^2} \sigma_\alpha\Big| + \Big|\dint \varepsilon[D^\alpha,u_{\rm E}]\nabla \sigma_{\rm E}\dfrac{\varepsilon^2T_{\rm app}}{(1+\varepsilon\sigma^\varepsilon)^2} \sigma_\alpha\Big|+ C\|U_E\|_{s,{\rm e}}^2\\
&\leq & \Big|\dint \dfrac{\varepsilon}{2}\nabla |\sigma_\alpha|^2\cdot u_{\rm E}\dfrac{\varepsilon^2T_{\rm app}}{(1+\varepsilon\sigma^\varepsilon)^2} \Big| + C\|U_E\|_{s,{\rm e}}^2\\
&= & \Big|\dint \dfrac{\varepsilon}{2} |\sigma_\alpha|^2 {\rm div}\Big( u_{\rm E}\dfrac{\varepsilon^2T_{\rm app}}{(1+\varepsilon\sigma^\varepsilon)^2} \Big) \Big| + C\|U_E\|_{s,{\rm e}}^2\leq C\|U_E\|_{s,{\rm e}}^2.
\end{eqnarray*}

According to the previous estimates, we arrive at
\begin{eqnarray}
\label{estcalR}\langle\mathcal{R}_\alpha,U_\alpha\rangle \leq C+C\|U_E\|_{s,{\rm e}}^2 +\eta\|\nabla u_\alpha\|_{L^2}^2+\eta \|\nabla T_0\|_{H^s}^2, \quad \mbox{ for small } \eta>0.
\end{eqnarray}

Finally, from \eqref{energyest}, \eqref{calL}, \eqref{singulartermend} and \eqref{estcalR}, we end up with the energy estimate
\begin{eqnarray*}
&&\dfrac{d}{dt}\dint (\frac{\varepsilon^2 T_{\rm app}}{(1+\varepsilon\sigma^\varepsilon)^2}\sigma_\alpha^2 +|u_\alpha|^2 +\dfrac{3}{2T_{\rm app}}T_\alpha^2+|\nabla\psi_\alpha|^2) \\
&&\qquad+\dint |u_\alpha|^2 +\dint |\nabla u_\alpha|^2 +\dint |\nabla T_\alpha|^2 \leq C (\|U_{\rm E}\|_{s,{\rm e}}^2 +\|\nabla\psi_0\|_{H^s}^2)+C.
\end{eqnarray*}
Taking the summation for all $0\leq |\alpha|\leq s$ and by using Gronwall's inequality, we have
\begin{eqnarray}
\label{energyestfinal} \dsup_{0\leq t\leq \tau} (\|U_{\rm E}\|_{s,{\rm e}}^2 +\|\nabla\psi_0\|_{H^s}^2) +\dint_0^\tau (\|u_{\rm E}\|_{H^{s+1}} +\|T_{\rm E}\|_{H^{s+1}}) \leq Ce^{C\tau}.
\end{eqnarray}

\section{Proof of the Main Result Theorem \ref{MT}}
By using the asymptotic expansion \eqref{expansionintro}, the Theorem \ref{thmincom}, and the  Lemma \ref{lem1f} we get the existence and uniqueness of classical solutions to the initial value problem for the system \eqref{scaledsystem}, \eqref{scaledsystemIC} and the solution satisfies
\begin{eqnarray*}
&\dsup_{0\leq t\leq \tau}\big[\|\varepsilon(\sigma^\varepsilon-\sigma_{\rm 1f})\|_{H^s}+\|u^\varepsilon-v-u_{\rm 1f}\|_{H^s} +\|T^\varepsilon -T\|_{H^s} +\|\nabla\psi^\varepsilon-\nabla\psi_{\rm 1f}\|_{H^s}\big] \leq C\varepsilon.&\\
&\|u^\varepsilon-v-u_{\rm 1f}\|_{L^2([0,\tau];H^{s+1})} +\|T^\varepsilon -T\|_{L^2([0,\tau];H^{s+1})} \leq C\varepsilon.
\end{eqnarray*}

\section{Appendix: Formal derivation of leading order oscillation}
We  follow the ideas of  Masmoudi's \cite{MasmoudiVP} and Schochet's \cite{Schochet} in order to deal with  the leading order oscillation system.
We start from \eqref{scaledsystem}, for simplicity, we drop all $\varepsilon$ in the superscripts.

First, it is reasonable to rewrite the mass conservation equation by using the electronic field
$$
\partial_t\nabla \psi +\dfrac{1}{\varepsilon}u= -(\sigma u)=-u\Delta\psi.
$$
Since there are no oscillation for the emperature equation, the leading order oscillations have to come from the equations from mass and velocity, namely
\begin{eqnarray}\label{eqnupsi}
\partial_t\left(\begin{array}{c} u\\ \nabla\psi\end{array}\right) +\dfrac{1}{\varepsilon}\left(\begin{array}{cc} 0&-\mathbb{I}\\ \mathbb{I}&0\end{array}\right) \left(\begin{array}{c} u\\ \nabla\psi\end{array}\right) = \left(\begin{array}{c} -u\cdot\nabla u +F(u)+G\\ -u\Delta \psi \end{array}\right)
\end{eqnarray}
where
\begin{eqnarray*}
F(u)&=&F(\varepsilon\sigma,T,u)=\dfrac{1}{1+\varepsilon\sigma}{\rm div}\Big[(1+\varepsilon\sigma) T(\nabla u+(\nabla u)^T-\dfrac{2}{3}{\rm div}u \mathbb{I})\Big]+u,\\
G&=&G(\varepsilon\sigma,T)=-\dfrac{\varepsilon T}{1+\varepsilon\sigma}\nabla \sigma -\nabla T.
\end{eqnarray*}

\subsection{Mapping $L$ and $e^{Lt}$}
Let $L$ be a linear mapping from $L^2(\T^d;\R^{2d})$ to itself, which is defined in the following way.

For any $v,e \in L^2(\T^d;\R^d)$ with ${\rm div}v={\rm div } e=0$, we set
\begin{eqnarray}
\label{defnL}& L\left(\begin{array}{c} v\\0\end{array}\right)=L\left(\begin{array}{c} 0\\ e\end{array}\right)=0,\quad&\\
\label{defnLnabla}&\mbox{ and }\quad L\left(\begin{array}{c} \nabla q\\ \nabla\phi\end{array}\right)
=\left(\begin{array}{cc} 0& -\I\\ \I & 0\end{array}\right)
\left(\begin{array}{c} \nabla q\\ \nabla \phi\end{array}\right).&
\end{eqnarray}

Since for any $ (u,E)^T\in L^2(\T^d;\R^{2d})$, it is well known that the Hodge decomposition yields to
$$
\left(\begin{array}{c} u\\ E\end{array}\right)=\left(\begin{array}{c} v\\ 0\end{array}\right)+\left(\begin{array}{c}0 \\ e\end{array}\right) + \left(\begin{array}{c} \nabla q\\ \nabla\phi\end{array}\right)
$$
where $v,e$ are divergence free and $q,\phi\in H^1(\T^d;\R^d)$,  $L$ is well defined by the identities \eqref{defnL}, \eqref{defnLnabla}.

The eigenvalues of $L$ are $0$, $i$ and $-i$, the corresponding eigenspaces are given by
\begin{eqnarray*}
E_0&=&\left\{\left(\begin{array}{c} v\\0\end{array}\right),\left(\begin{array}{c} 0\\e\end{array}\right), \mbox{ where } v,e \in L^2(\T^d;\R^d) \mbox{ and } {\rm div} v={\rm div}e =0\right\}\\
E_i &=& \left\{\left(\begin{array}{c} \nabla q\\ -i\nabla q\end{array}\right), \mbox{ where } q\in H^1(\T^d;\R^d)\right\}\\
E_{-i} &=& \left\{\left(\begin{array}{c} \nabla q\\ i\nabla q\end{array}\right), \mbox{ where } q\in H^1(\T^d;\R^d)\right\}.
\end{eqnarray*}

For any $(u,E)^T\in L^2(\T^d;\R^d)$ with $u=v+\nabla q$, $E=e+\nabla\phi$, the projection operators according to $L$ are given by
\begin{eqnarray*}
P_0\left(\begin{array}{c} u\\E\end{array}\right)&=& \left(\begin{array}{c} v\\e\end{array}\right),\\
P_i\left(\begin{array}{c} u\\E\end{array}\right)&=& \dfrac{1}{2}\left(\begin{array}{c} \nabla q+ i\nabla\phi\\-i\nabla q+\nabla\phi\end{array}\right),\\
P_{-i}\left(\begin{array}{c} u\\E\end{array}\right)&=& \dfrac{1}{2}\left(\begin{array}{c} \nabla q- i\nabla\phi\\i\nabla q+\nabla\phi\end{array}\right).
\end{eqnarray*}

Then $(I-P_0)e^{\tau L}$ is defined by
$$
(I-P_0)e^{\tau L}\left(\begin{array}{c} u\\E\end{array}\right)=\left(\begin{array}{c}\cos \tau \nabla q-\sin\tau\nabla\phi \\ \cos\tau\nabla\phi+\sin\tau\nabla q\end{array}\right)=\left(\begin{array}{cc}\cos \tau  & -\sin\tau \\ \sin\tau &\cos\tau\end{array}\right) \left(\begin{array}{c}\nabla q \\ \nabla\phi\end{array}\right).
$$

\subsection{Asymptotic expansion by fast and slow time}

Assume that we have the expansion of solution, let $\tau=\frac{t}{\varepsilon}$
\begin{eqnarray*}
u(x,t) &=& u^{(0)}(x,t,\tau) +\varepsilon u^{(1)}(x,t,\tau) +\varepsilon^2 u^{(2)}(x,t,\tau)+\cdots,\\
\nabla \psi(x,t) &=& \nabla\psi^{(0)}(x,t,\tau) +\varepsilon \nabla\psi^{(1)}(x,t,\tau) +\varepsilon^2 \nabla\psi^{(2)}(x,t,\tau)+\cdots,\\
T(x,t) &=& T^{(0)}(x,t) +\varepsilon T^{(1)}(x,t,\tau) +\varepsilon^2 T^{(2)}(x,t,\tau)+\cdots,\\
\sigma(x,t) = \Delta \psi(x,t) &=& \Delta\psi^{(0)}(x,t,\tau) +\varepsilon \Delta\psi^{(1)}(x,t,\tau)+\varepsilon^2 \Delta\psi^{(2)}(x,t,\tau)+\cdots,
\end{eqnarray*}
By putting this ansatz into our system and comparing the orders of $\varepsilon$, we have the following systems.

The $O(\frac{1}{\varepsilon})$ order terms give
\begin{eqnarray*}
\partial_\tau u^{(0)} -\nabla\psi^{(0)} =0,\\
\partial_\tau\nabla\psi^{(0)} +\nabla u^{(0)}=0.
\end{eqnarray*}
The $O(1)$ order terms give
\begin{eqnarray*}
\partial_t u^{(0)} +\partial_\tau u^{(1)} -\nabla\psi^{(1)}& =& -u^{(0)}\cdot\nabla u^{(0)} + F(0,u^{(0)}, T^{(0)}) +G(0,T^{(0)})\\
\partial_t \nabla \psi^{(0)} +\partial_\tau\nabla\psi^{(1)}+u^{(1)} &=& -u^{(0)}\Delta\psi^{(0)}.
\end{eqnarray*}
By using operator $L$ we defined before, we are able to rewrite this system into
\begin{eqnarray*}
&&\left(\begin{array}{c}u^{(1)} \\ \nabla\psi^{(1)}\end{array}\right) = e^{-L\tau}\left(\begin{array}{c}u^{(1)}(x,t,0) \\ \nabla\psi^{(1)}(x,t,0) \end{array}\right) \\
&&\qquad+e^{-L\tau}\dint^\tau_0 e^{Ls}\left(\begin{array}{c} -\partial_tu^{(0)}-u^{(0)}\cdot\nabla u^{(0)} + F(0,u^{(0)}, T^{(0)}) +G(0,T^{(0)}) \\ -\partial_t\nabla\psi^{(0)}-u^{(0)}\Delta\psi^{(0)}\end{array}\right)ds.
\end{eqnarray*}
To have the expansion meaningful, the first order term should be $u^{(1)},\nabla\psi^{(1)}\ll \tau$ as $\tau\rightarrow\infty$. Since it is obvious that the initial data term
\begin{eqnarray*}
\dfrac{1}{\tau}e^{-L\tau}\left(\begin{array}{c}u^{(1)}(x,t,0) \\ \nabla\psi^{(1)}(x,t,0) \end{array}\right) \rightarrow 0,
\end{eqnarray*}
then it is necessary to check that
\begin{eqnarray}
\label{osc1iff}
\dfrac{1}{\tau}\dint^\tau_0 e^{Ls}\left(\begin{array}{c} -\partial_tu^{(0)}-u^{(0)}\cdot\nabla u^{(0)} + F(0,u^{(0)}, T^{(0)}) +G(0,T^{(0)}) \\ -\partial_t\nabla\psi^{(0)}-u^{(0)}\Delta\psi^{(0)}\end{array}\right)ds \rightarrow 0, \,\mbox{ as }\tau\rightarrow 0.
\end{eqnarray}

From the eigenvalue analysis of operator $L$, we know that for any given $(u^{(0)},\nabla\psi^{(0)})$, there exists $(\nabla q,\nabla\phi)$ such that
$$
P_0\left[e^{\frac{t}{\varepsilon}L}\left(\begin{array}{c}u^{(0)}\\ \nabla\psi^{(0)}\end{array}\right)\right]=\left(\begin{array}{c}v\\ 0\end{array}\right) \quad \mbox{ and }\quad (I-P_0)\left[e^{\frac{t}{\varepsilon}L}\left(\begin{array}{c}u^{(0)}\\ \nabla\psi^{(0)}\end{array}\right)\right]=\left(\begin{array}{c}\nabla q\\ \nabla\phi\end{array}\right),
$$
In other words we have
\begin{eqnarray*}
\left(\begin{array}{c}u^{(0)}\\ \nabla\psi^{(0)}\end{array}\right) = \left(\begin{array}{c}v\\ 0\end{array}\right) +e^{-\frac{t}{\varepsilon}L}\left(\begin{array}{c}\nabla q\\ \nabla\phi\end{array}\right)
%=\left(\begin{array}{c}v+ \cos(-\frac{t}{\varepsilon})\nabla q -\sin(-\frac{t}{\varepsilon})\nabla\phi\\ \cos(-\frac{t}{\varepsilon})\nabla \phi %+\sin(-\frac{t}{\varepsilon})\nabla q\end{array}\right)
=\left(\begin{array}{c}v+ \cos\frac{t}{\varepsilon}\nabla q +\sin\frac{t}{\varepsilon}\nabla\phi\\ \cos\frac{t}{\varepsilon}\nabla \phi -\sin\frac{t}{\varepsilon}\nabla q\end{array}\right).
\end{eqnarray*}

By applying operator $P_0$ on both sides of \eqref{osc1iff}, we will have the incompressible limit.

Applying operator $I-P_0$ on both sides of \eqref{osc1iff}, we will get our leading order fast oscillation system, which will be shown in the following calculation.

We will calculate the following limit
\begin{eqnarray*}
&&\dlim_{s\rightarrow +\infty}\dfrac{1}{s}\dint^s_0 d\tau (I-P_0) e^{\tau L}\left(\begin{array}{c} QI_1\\ QI_2 \end{array}\right)\\
 &=&\dlim_{s\rightarrow +\infty}\dfrac{1}{s}\dint^s_0d\tau \left(\begin{array}{c} \cos\tau QI_1 -\sin\tau QI_2\\ \cos\tau QI_2 +\sin \tau QI_1\end{array}\right)\\
&=&\dfrac{1}{2\pi}\dint^{2\pi}_0d\tau \left(\begin{array}{c} \cos\tau QI_1 -\sin\tau QI_2\\ \cos\tau QI_2 +\sin \tau QI_1\end{array}\right),
\end{eqnarray*}
where
\begin{eqnarray*}
I_1&=&-(v+ \cos\tau\nabla q +\sin\tau\nabla\phi)\cdot\nabla (v+ \cos\tau\nabla q +\sin\tau\nabla\phi) \\
&& \hspace{2cm} +F(v+ \cos\tau\nabla q +\sin\tau\nabla\phi)+G\\
I_2&=&-(v+ \cos\tau\nabla q +\sin\tau\nabla\phi)\cdot (\cos\tau\Delta \phi -\sin\tau\Delta q).
\end{eqnarray*}

Thus the desired oscillation equations are
\begin{eqnarray}
\label{oscq}
2\partial_t\nabla q &= &Q\left\{-(\nabla q\cdot\nabla)v -(v\cdot\nabla)\nabla q + v\Delta q +{\rm div} [T \mathbb{S}(\nabla q)]\right\}-\nabla q,\\
\label{oscphi}
2\partial_t\nabla \phi &= &Q\left\{-(\nabla \phi\cdot\nabla)v -(v\cdot\nabla)\nabla \phi + v\Delta \phi+ {\rm div} [T \mathbb{S}(\nabla \phi)]\right\}-\nabla \phi.
\end{eqnarray}
with initial data
\begin{eqnarray}
\label{oscIC}\nabla q |_{t=0}=Q u_{\rm I}(x),\quad \nabla \phi|_{t=0}=\nabla\psi_{\rm I}(x),
\end{eqnarray}
where $\psi_{\rm I}(x)$ satisfies $\Delta \psi_{\rm I}(x)=\sigma_{\rm I}(x)$.

It is straight forward to prove the following result by using energy estimates.
\begin{lem}
\label{lemexistenceqphi} For any fixed $\tau>0$, $\forall s>\frac{d}{2}+1$.
Given $v,T\in L^\infty([0,\tau];H^s)$, $T\geq T_{\rm L}>0$ and $Q u_{\rm I}(x),\nabla\psi_{\rm I}(x)\in H^s$.
Problem \eqref{oscq}\eqref{oscphi}\eqref{oscIC} has a unique solution $(\nabla q, \nabla\phi)$ such that
\begin{eqnarray}
\dsup_{0\leq t\leq\tau}\|\nabla q(\cdot, t),\nabla \phi(\cdot, t)\|_{H^s} +\|\nabla q,\nabla \phi\|_{L^2([0,\tau];H^{s+1})}\leq C(\tau) \|Q u_{\rm I},\nabla \psi_{\rm I}\|_{H^s}.
\end{eqnarray}
\end{lem}

\end{document}